\documentclass{article}

\usepackage{amscd}

\usepackage{amssymb,latexsym}

\begin{document}

\pagestyle{myheadings}

\markright{MARCELO GOMEZ MORTEO}

\title  {The Farrell-Jones Isomorphism Conjecture in K-Theory }
\author{ Marcelo Gomez Morteo }

\maketitle
\vspace{16pt}

\begin{abstract}

\vspace{16pt}

We prove that the Farrell-Jones isomorphism conjecture for non-connective algebraic K-theory for a discrete group $G$ and a coefficient ring $R$ holds true if $G$ belongs to the class of groups acting on trees, under certain conditions on $G$ (see theorem 0.5 below ) and if the coefficient ring $R$ is either regular or hereditary, depending on the structure of $G$. Our result is weaker than the result that has been established in [15] which says that these groups, verify the conjecture for any coefficient ring, see remark 0.6 below.

\end{abstract}

\vspace{30pt}

\textbf{0 Introduction:}

\vspace{30pt}

 Let $G$ denote a discrete group which is vitrtually torsion free. A virtually torsion free group $G$ is a group $G$ containing a subgroup which is torsion free and of finite index. There are many examples of such class of groups: for example any arithmetic group ( see [9] page 108) or any Fuchsian group (See [7] page 100) belongs to this class. A finitely generated Fuchsian group has a presentation given by ([8]),

\vspace{20pt}

\[<a_{1},a_{2},....a_{g},b_{1},b_{2},...b_{g},c_{1},c_{2},...c_{n}/ c_{1}^{e_{1}}=c_{2}^{e_{2}}=....=c_{n}^{e_{n}}=1,c_{1}c_{2}...c_{n}\Pi_{i=1}^{g}[a_{i},b_{i}]=1>\]

\vspace{20pt}

where $[a_{i},b_{i}]=a_{i}b_{i}a_{i}^{-1}b_{i}^{-1}$. A torsion free subgroup of a Fuchsian group is either a free group or the fundamental group

\vspace{20pt}

\[<a_{1},a_{2},....a_{g},b_{1},b_{2},...b_{g}/ \Pi_{i=1}^{g}[a_{i},b_{i}]=1>\]

\vspace{20pt}

Given a ring with unit $R$, we say that $G$ has finite $R$-cohomological dimension, written $cd_{R}(G) < \infty$ if there exists a projective resolution of $R$ by right projective $R[G]$-modules with finite projective dimension. If $R=\mathbb{Z}$  the ring of integers, then we note $cd_{R}(G)=cd(G)$.
The virtual $R$-cohomological dimension for $R$ commutative is defined as $vcd_{R}(G)=cd_{R}(\Delta)$ for any torsion subgroup $\Delta \subset G$ with finite index $[G:\Delta]$. The virual cohomological dimension is well defined since for any torsion free subgroups $\Delta$ and $\Gamma$ of $G$ their $R$-cohomological dimensions are coincident by a theorem of Serre (see [1] page 30).We note $vcd(G)$ for the virtual cohomological dimension of $G$ associated to the ring $\mathbb{Z}$. A family of subgroups of $G$ is a set of (closed) subgroups of $G$ which is closed under conjugation and finite intersections. The examples which we will be working with are the family $\cal{FIN}$ of finite subgroups of $G$ and the family $\cal{VCYC}$ of virtually cyclic subgroups of $G$. A subgroup $V$ is called virtually cyclic if it is either finite or if it contains an infinite cyclic subgroup of finite index. Let $E_{\cal{FIN}}(G)$ be a model for the family $\cal{FIN}$ (see [6] page 289). If $G$ is virtually torsion free, then $vcd(G)< \infty$ if and only if there exists a finite dimensional contractible simplicial complex $X$ on which $G$ acts properly (and simplicially). Moreover, the space $X$ can be chosen so that the fixed point set $X^{H}$ is contractible for every finite subgroup $H\subset G$ (See [1], proposition 2.1).This means ( see [6] page 290) that $X$ is a finite model for $G$.

\vspace{10pt}

 Arithmetic groups have this property. A subgroup $\Gamma$ of $G(\mathbb{Q})$ is arithmetic if it is commensurable with $G(\mathbb{Z})$, that is, if $ \Gamma \cap G(\mathbb{Z})$ has finite index both in $\Gamma$ and in $G(\mathbb{Z}$. A group is arithmetic if it can be embedded as an arithmetic subgroup of $G(\mathbb{Q})$ for some $\mathbb{Q}$-algebraic subgroup $G(\mathbb{Q})$ of $GL_{n}$. Then any subgroup of finite index in $\Gamma$ is also an arithmetic group. An arithmetic subgroup $\Gamma$ in a semisimple connected linear $\mathbb{Q}$-algebraic group has a finite model for $E_{\cal{FIN}}(\Gamma)$. Let $G(\mathbb{R})$ be the $\mathbb{R}$-points of a semisimple $\mathbb{Q}$-group $G(\mathbb{Q})$, and let $K \subset G(\mathbb{R})$ be a maximal compact group. If $\Gamma \subset G(\mathbb{Q})$ is an arithmetic group, then $G(\mathbb{R})/K$ with the left $\Gamma$-action is a model for $E_{\cal{FIN}}(\Gamma)$. However it is not necessarily a finite model. There is a finite model though given by the Borel-Serre compactification of $G(\mathbb{R})/K$,( see [6]). Since as we have already mentioned, any arithmetic groups is virtually torsion free, then by proposition 2.1 in [1] it must be of finite virtual cohomological dimension. In any case we have the following result stated in [6] page 296.

 \vspace{20pt}

 \textbf{Theorem 0.1}:\emph{Let $L$ be a Lie group with finitely many path components. Let $K$ be a maximal compact subgroup. Let $G \subset L$ be a discrete subgroup of $L$. Then $L/K$ with the left $G$-action is a model for $E_{\cal{FIN}}(G)$. Suppose additionally that $G$ is virtually torsion free, then we have for its virtual cohomological dimension,
 }
 \vspace{20pt}

 \[vcd(G)\leq dim(L/K)\]

 \vspace{20pt}

 \emph{and equality holds if and only $G\setminus L$ is compact.
 }
\vspace{20pt}

For a maximal compact subgroup $K \subset L$ it is known that the homogeneous space $X=L/K$ is diffeomorphic to the Euclidean space $\mathbb{R}^{d}$ where $d=dim L-dim K$ and that $G$ acts properly on $X$. In particular every isotropy group $G_{x}$ is finite.Moreover $G\setminus L$ is compact if and only $G\setminus X$ is compact,(see [1] page 28 and also [9] page 111).

\vspace{15pt}

Since any Fuchsian group $F$ is virtually torsion free and is included in the Lie group $PSL(2,\mathbb{R})$, then, $vcd(F)\leq 2$. In this case $X=\mathbb{H}$ where $\mathbb{H}$ is the hyperbolic plane and by the argument of the last paragraph, $vcd(F)\leq 1$ if and only if $F$ is not cocompact, ie if $F\setminus \mathbb{H}$ is not compact. For an arithmetic Fuchsian group this is the case if it is derived from a quaternion algebra isomorphic to the 2x2-matrices with coefficients in $\mathbb{Q}$ noted $M(2,\mathbb{Q})$, (see [10] page 113)

\vspace{15pt}

The virtually torsion free groups $G$ with $vcd(G)\leq 1$ carry special information: such a $G$ is the fundamental group of a connected graph of finite groups of bounded order, (see [1], page 32 and [5]).By [5], page 119, this fact is equivalent to the statement that $G$ has a free subgroup of finite index. In particular a non-cocompact Fuchsian group has a free subgroup of finite index. A discrete group $G$ which is virtually torsion free has virtual cohomological dimension bounded by the dimension of any model for that group with respect to the family of finite subgroups, (see [6] page 297). Therefore if $G$ acts on a tree and is virtually torsion free it has $vcd \leq 1$ and has therefore a free subgroup of finite index, ie $G$ is virtually free. conjecture.We say that a group $G$ verifies property FA if every isometric action of $G$ on a tree fixes a vertex or an edge, (see [16]) The point is that by the structure theorem 0.4 (see below ) and by Lemma 7.4 in [11] page 192, any infinite discrete group $G$ acting properly on a tree $T$ (ie, does not verify property FA), has a free subgroup of finite index, has therefore virtual cohomological dimension bounded by 1, and is virtually torsion free. Moreover if $G$ is finitely generated, then it is an hyperbolic group and therefore verifies the Farrell-Jones isomorphism conjecture. Finitely generated groups which are virtually free are hyperbolic (see [6] page 312). Also finitely generated Fuchsian groups are also hyperbolic, (see [2] page 450) and then they verify the Farrell-Jones isomorphism. 

\vspace{10pt}

Suppose $G$ is a discrete infinite group which is not locally finite. If $G$ is finite or locally finite, it satisfies the Farrell-Jones isomorphism conjecture. Then adapting the theorem 5.1 on page 83 [5] we obtain:

\vspace{20pt}

\textbf{Theorem 0.4:} \emph{Suppose $G$ is an infinite discrete group. The following two statements are equivalent:}

\vspace{15pt}

\emph{a) $G$ acts properly on a tree $T$ and hence does not verify property FA.
}
\vspace{15pt}

\emph{b) $G$ is the amalgamated product $G=A*_{H}B$ with $A,B,H$ finite groups, or}

\vspace{10pt}

\emph{$G$ is an HNN extension $G=<\alpha,\beta:H \mapsto A, t>$ with $A,H$ finite groups.}

\vspace{10pt}

\emph{c) $T/G$ has precisely one edge.}

\vspace{20pt}

Remark 0.6: a)An infinite discrete group $G$ acting properly on a tree (see [2], [16] for the definition) also acts cocompactly on that tree by c) of theorem 0.4. A tree is a one dimensional Cat(0) space (see [2]), therefore since the Farrell-Jones isomorphism conjecture has been established for $G$-Cat(0) groups, ie groups acting properly and cocompactly on a finite dimensional Cat(0) space ( see [15])it follows that for groups $G$ under the conditions of theorem 0.4 verify the Farrell-Jones conjecture for any coefficient ring.

\vspace{10pt}

b) Some of the groups not enjoying the hypothesis of theorem 0.4, ie having the FA property are the finitely generated torsion groups or the groups $SL(n,\mathbb{Z})$ for all $n \geq 3$.

\vspace{20pt}

We are interested in the following isomorphism conjectures:

\vspace{10pt}

\textbf{Conjecture 0.2}: \emph{($K$-theoretic Farrell-Jones conjecture for torsion free groups). The $K$-theoretic Farrell-Jones conjecture with coefficients in the (right )regular ring $R$ for the torsion free group $G$ predicts that the assembly map
}
\vspace{20pt}

\[H_{n}(BG,K(R)) \mapsto K_{n}(R[G])\]

\vspace{15pt}

\emph{is bijective for all $n \in \mathbb{Z}$.}

\vspace{10pt}

 Remember that BG is the classifying space of the group $G$, $K(R)$ is the non-connective algebraic $K$-theory spectrum of the ring $R$, $K_{n}(R[G])$ is the algebraic $K$-theory of the group ring $R[G]$. $R$ being (right )regular means that $R$ is (right) noetherian and and every finitely generated (right) $R$-module possesses a finite projective resolution, (see [12] page 122).

\vspace{20pt}

\textbf{Conjecture 0.3:} ($K$-theoretic Farrell-Jones conjecture). \emph{The $K$-theoretic Farrell-Jones conjecture with coefficients in a ring with unit $R$ for the group $G$ predicts that the assembly map}

\vspace{20pt}

\[H_{n}^{G}(E_{\cal{VCYC}}(G),K_{R}) \mapsto K_{n}(R[G])\]

\vspace{15pt}

\emph{which is the map induced by the projection $E_{\cal{VCYC}}(G) \mapsto * $ is bijective for all $n \in \mathbb{Z}$}.

\vspace{15pt}

Here $H_{*}^{G}(-,K_{R})$ is the $G$-equivariant homology associated to a certain functor $K_{R}:\cal{GROUPOIDS} \mapsto \cal{SPECTRA}$ which satisfy for every group $G$ and all $n\in \mathbb{Z}$, $\pi_{n}(K_{R}(G) \simeq K_{n}(R[G])$, (see [6] pages 299-300).

\vspace{10pt}

These conjectures are an important tool for the computation of $K$-groups. Computation of $K$-groups is a difficult task, but with the disposal of the isomorphism conjectures they can be accomplished through the computation of an homology which is usually easier to calculate with the aid of adequate spectral sequences.

\vspace{10pt}

In section 1 we give a partial result on these conjectures:

\vspace{20pt}

\textbf{Theorem 0.5:} \emph{Let $G$ be an infinite discrete group acting properly on a tree. With the notation of theorem 0.4 , if:}

\vspace{10pt}

\emph{a)$R$ is (right) hereditary ring and $H$ is non trivial and is without $R$-torsion,
}
\vspace{10pt}

\emph{ b) $R$ is (right)regular and $H$ is the trivial group,}

\vspace{15pt}

\emph{then the Farrell-Jones conjecture is verified for $G$ and $R$.}

\vspace{20pt}

A group is without $R$ torsion, and sometimes it is then called an $R^{-1}$-group, if the orders of all finite subgroups of $G$ are invertible. A ring is (right) hereditary if all its (right) ideals are (right) projective.

\vspace{10pt}

A ring $R$ is (right) coherent if every (right) finitely presentable $R$-module has a projective resolution by finitely generated projective modules, and (right )regular coherent if this resolution can be taken to be finite dimensional. Therefore a ring $R$ is (right) regular coherent if every finitely presentable right $R$-module has a finite dimensional projective resolution by projective $R$-modules which are finitely generated,(See [13] page 160). Observe that with the above definition of a (right) regular ring, a ring is (right) regular if and only if it is (right) regular coherent and (right) noetherian. With the hypothesis of theorem 0.4, $R[H]$ is (right) regular coherent.We use to prove that $R[H]$ is (right) regular coherent a result proven [14] theorem 1, which states that if $H$ is a finite group without $R$- torsion and $R$ is (right) hereditary, then $R[H]$ is (right) hereditary. A (right) hereditary ring is always (right) regular coherent. Incidentally, the fact proven in theorem 1 in [14] mentioned in this paragraph is a generalization of another more classical theorem. That theorem affirms that if $R$ is a Dedekind domain (which is the counterpart in the commutative setting of an hereditary ring )the group ring $R[G]$ is a maximal order if and only if the order of $G$ is a unit in $R$. The point is that a maximal order in a semisimple algebra is both right and left hereditary, ie is hereditary. The converse is false, ( [4], vol 1, pages 564,565 ).

\vspace{20pt}

An immediate consequence of this theorem is that if $G$ is a free group then the Farrell-Jones isomorphism conjecture is verified for all(right)regular rings, (see [5] page 119). In particular since the Farrell-Jones isomorphism conjecture and the Farrell-Jones conjecture are equivalent for torsion free groups, (see [6]), then both isomorphism conjectures are verified in this special case.

\vspace{15pt}

Another consequence of this theorem is that the Farrell-Jones isomorphism conjecture is true for finitely generated non-cocompact Fuchsian groups with regular coefficient rings, since these groups can be written as free products on a finite number of finite cyclic groups plus a free group, (see [8]) and also for  amalgamated finite groups and HNN extensions of finite groups if the group is without $H$-torsion (with notation of theorem 0.4) and the coefficient ring $R$ is (right) hereditary. See the examples in section 1.

\vspace{30pt}

\textbf{1 Proof of Theorem 0.5 and its consequences}

\vspace{30pt}

\textbf{Definition 1.1:} \emph{Following [3] definition 4.1, for any discrete group $G$ acting on a tree $T$ we consider the projection $pr:T \mapsto *$ and we say that the homology theory $H_{*}^{?}(-,K_{R})$ has the tree property if}

\vspace{20pt}

\[H_{n}^{G}(T,K_{R}) \mapsto H_{n}^{G}(*,K_{R})\]

\vspace{20pt}

\emph{is an isomorphism for all $n \in \mathbb{Z}$.}

\vspace{20pt}

An infinite discrete group $G$ acting on a tree is either an amalgamated product or an HNN extension by theorem 0.4 . By lemma 5.1 in [3] $H_{*}^{?}(-,K_{R})$ has the tree property if and only conditions a) and b) stated below are satisfied:

\vspace{15pt}

 a) When $G$ is an amalgamated product, ie a pushout

\vspace{20pt}

\[
\begin{CD}
 H @>>> G_{1}\\
@VVV @VVV\\
 G_{2} @>>> G
\end{CD}
\]

\vspace{20pt}

this pushout induces an homotopy cocartesian diagram noted $D_{1}$

\vspace{20pt}

\[
\begin{CD}
 K_{R}(H)\vee K_{R}(H) @>>> K_{R}(G_{1})\vee K_{R}(G_{2})\\
@VVV @VVV\\
 K_{R}(j):K_{R}(H) @>>> K_{R}(G)
\end{CD}
\]

\vspace{20pt}

where $j:H \mapsto G$ is defined to be $j_{1}\circ i_{1}=j_{2}\circ i_{2}$, where $i_{l}:H \mapsto G_{l}$, $l=1,2$ and $j_{l}:G_{l} \mapsto G$, $l=1,2$ are the injective maps of the pushout square which define the amalgamated product.

\vspace{15pt}

b) If $G$ is an HNN extension, this extension having injective morphisms $s_{1},s_{2}:H \mapsto L$  associated to $G$, and an inclusion $i:L \mapsto G$, these maps induce an homotopy cocartesian diagram noted $D_{2}$

\vspace{20pt}

\[
\begin{CD}
 K_{R}(H)\vee K_{R}(H) @>>> K_{R}(L)\\
@VVV @VVV\\
 K_{R}(j):K_{R}(H) @>>> K_{R}(G)
\end{CD}
\]

\vspace{15pt}

where $j=i\circ s_{1}=i\circ s_{2}$ See remark 5.6 in [3]

\vspace{20pt}

It turns out that Waldhausen in [13] defines certain spectra, called Nil-Spectra, noted $Nil^{-\infty}(R,X,Y)$,(see [3], def 9.4), where $R$ is a ring and $X$ and $Y$ are $R$-bimodules. With the notation given above, the injective morphisms of groups $i_{1},s_{1},i_{2},s_{2}$ induce the respective ring morphisms $R(i_{1}),R(s_{1}),R(i_{2}),R(s_{2})$ on group rings. This ring morphisms are pure and free, (see [3] below definition 9.4). An inclusion $\alpha:C \mapsto A$ of rings is called pure if $A=\alpha(C)\oplus A^{\bullet}$ and $A^{\bullet}$ is a $C$-bimodule, and moreover it is pure and free if $A^{\bullet}$ is free. In our setting we therefore have that $R[G_{l}]=R(i_{l})(H)\oplus R[G_{l}]^{\bullet}$ for $l=1,2$ and $R[L]=R(s_{l})(H)\oplus R[L]^{\bullet}$ for $l=1,2$ with $R[G_{1}]^{\bullet},R[G_{2}]^{\bullet},R[L]^{\bullet}$ free $R$-bimodules.

\vspace{10pt}

The left upper corner of the diagram $D_{1}$ defined above if replaced by

\vspace{10pt}

\[Nil^{-\infty}(R[H],R[G_{1}]^{\bullet},R[G_{2}]^{\bullet})\]

\vspace{10pt}

gives an homotopy cartesian diagram, (see [3] theorem 10.2).

\vspace{20pt}

Similarly in [13] a second Nil-Spectra noted $Nil^{-\infty}(R,X,Y,Z,W)$ is defined, where $R$ is a ring and $X,Y,Z,W$ are $R$ bimodules. Let $\alpha,\beta:C \mapsto A$ be pure and free. The Laurent extension with respect to $\alpha,\beta$ is the universal ring $R=_{\alpha}A_{\beta}(t,t^{-1})$ that contains $A=_{\alpha}A_{\beta}$ with this notation justified by the left and right actions on $A$ given by $\alpha,\beta$ respectively. Here $t$ is an invertible element and satisfies for all $c \in C$

\vspace{15pt}

\[\alpha(c)t=t\beta(c)\]

Its existence is explained in [13], page 149. In our setting $\alpha=R(s_{1})$ and $\beta=R(s_{2})$ define a Laurent extension with $C=R[H]$ and $A=R[L]$. This means that applying $R$ to an HNN extension gives a Laurent extension. Write $R[L]=R(s_{1})(H)\oplus R[L]^{\bullet}$ and $R[L]=R(s_{2})(H)\oplus R[L]^{\bullet\bullet}$.

\vspace{10pt}

The left upper corner of the diagram $D_{2}$, if replaced by

\vspace{20pt}

\[Nil=Nil^{-\infty}(R[H],_{R(s_{1})}R[L]_{R(s_{1})}^{\bullet},_{R(s_{2})}R[L]_{R(s_{2})}^{\bullet\bullet},
_{R(s_{2})}R[L]_{R(s_{1})},_{R(s_{1})}R[L]_{R(s_{2})})\]

\vspace{15pt}

gives an homotopy cartesian diagram. If $K^{-\infty}$ is the standard notation for non-connective $K$-theory, we have weak equivalence of spectra between $K^{-\infty}(R[H])$ and $K_{R}(H)$ for any group $H$. Also, by remarks 10.3 and 10.7 in [3], there are weak equivalences of spectra given by

\vspace{20pt}

\[K^{-\infty}(R[H])\vee K^{-\infty}(R[H]) \mapsto Nil^{-\infty}(R[H],R[G_{1}]^{\bullet},R[G_{2}]^{\bullet})\]

\vspace{20pt}

\[K^{-\infty}(R[H])\vee K^{-\infty}(R[H]) \mapsto Nil\]

\vspace{15pt}

if $R[H]$ is a ( right)regular coherent ring.

\vspace{10pt}

Therefore if $R[H]$ is (right) regular coherent, it implies that a) and b) are satisfied since in the first place these conditions are satisfied by replacing the left upper corners with the respective Nil Spectra, and moreover we know that we can replace the Nil Spectra by $K^{-\infty}(R[H])\vee K^{-\infty}(R[H])$ because the weak equivalences just stated hold if $R[H]$ is regular coherent. We can now prove theorem 0.5:

\vspace{20pt}

 By [14] theorem 1, and the hypothesis of theorem 0.5, $R[H]$ is (right) regular coherent. Then conditions a) and b) hold. If these conditions hold then the homology theory $H_{*}^{?}$ has the tree property. Since $G$ has as a model for the family $\cal{FIN}$ a tree $T$, we obtain for all $n \in \mathbb{Z}$,

\vspace{20pt}

\[H_{n}^{G}(E_{\cal{FIN}}(G),K_{R}) \mapsto K(R[G])\]

\vspace{20pt}

 Now for every virtually finite subgroup $V$ of $G$ either $V$ is finite, or $V$ is infinite and discrete and acts on a tree, (see [11], theorem 3.7 page 157). In any case we get for all $V \in \cal{VCYC}$,

 \vspace{20pt}

 \[H_{n}^{V}(E_{\cal{FIN}\cap V}(V),K_{R}) \mapsto K(R[V])\]

 \vspace{15pt}

 for all $n \in \mathbb{Z}$. Therefore by the transitivity principle stated in [6] page 309, we derive that for all $n \in \mathbb{Z}$,

 \vspace{20pt}

 \[H_{n}^{G}(E_{\cal{FIN}}(G),K_{R}) \mapsto H_{n}^{G}(E_{\cal{VCYC}}(G,K_{R})\]

 \vspace{15pt}

 is an isomorphism. Hence the Farrell-Jones isomorphism conjecture is established and theorem 0.5 is proven.

 \vspace{40pt}

 Example 1: The group $SL(2,\mathbb{Z})$ which is equal to $\mathbb{Z}/4\mathbb{Z}*_{\mathbb{Z}/2\mathbb{Z}}\mathbb{Z}/6\mathbb{Z}$. If we select a Dedekind domain $R$ then it is hereditary. If moreover 2 is invertible in this Dedekind ring, the Farrell-Jones conjecture holds true. Of course this is not a very strong statement since $SL(2,\mathbb{Z})$ is hyperbolic ([6] page 312) and consequently the Farrell-Jones conjecture is verified for any coefficient ring on this group.

 \vspace{20pt}

 Example 2: The group $PSL(2,\mathbb{Z})$ which is the free product $\mathbb{Z}/2\mathbb{Z}*\mathbb{Z}/3\mathbb{Z}$ verifies the Farrell-Jones conjecture for (right) regular coefficient rings since here the amalgamation group is the trivial group. This is a Fuchsian group with non compact fundamental region and is therefore non- cocompact, (see [10]). Observe again that this is not a strong statement since finitely generated Fuchsian groups are hyperbolic and therefore verify the Farrell-Jones isomorphism conjecture for any coefficient ring. As mentioned above with theorem 0.5 we prove that any non-cocompact finitely generated Fuchsian group verifies the Farrell-Jones isomorphism conjecture for any (right) regular ring.
 \vspace{20pt}

 Remark: Moreover, in example 1 or in example 2, we cannot affirm that having an isomorphism conjecture for the family $\cal{FIN}$ we have a stronger isomorphism conjecture than for the family $\cal{VCYC}$ since as we have seen, the transitivity principle applies.

 \vspace{20pt}

 Example 3: Any amalgamated product $G$ with finite amalgamation group $H$ or any HHN extension $G$ which extends through a finite group $H$ satisfies the Farrell'Jones isomorphism conjecture for coefficient rings wich are (right) hereditary and if $H$ is without $R$-torsion. If finitely generated, this groups always have a finite rank free subgroup of finite index,([11], theorem 7.4 page 192) and therefore are hyperbolic groups and verify the Farrell-Jones isomorphism conjecture. Still, $G$ could be non-finitely generated, but in this case our statement is still not as strong as the general result by remark 0.6.

\vspace{16pt}

\emph{E-mail address}: valmont8ar@hotmail.com

\end{document}